\documentclass[12pt]{article}
\usepackage{latexsym, amssymb, amsmath, amscd, amsfonts, epsfig, graphicx, colordvi,verbatim,ifpdf,tikz,extarrows}
\usepackage{amsfonts, amsmath, amssymb}
\usepackage{amssymb,amsfonts,amsmath,latexsym,epsfig,cite, psfrag,eepic,color}
\usepackage{amscd,graphics}
\usepackage{latexsym, amssymb,  amsmath,amscd, amsfonts, epsfig, graphicx, colordvi,amsthm}
\usepackage{ytableau}

\usepackage[colorlinks,linkcolor=blue,anchorcolor=blue,citecolor=blue]{hyperref}

\usepackage{graphicx}
\usepackage{color}
\usepackage{ifpdf}
\usepackage{fancybox}

\usepackage{float}

\newtheorem{thm}{Theorem}[section]

\newtheorem{core}[thm]{Corollary}

\def\qed{\nopagebreak\hfill{\rule{4pt}{7pt}}
\medbreak}

\hypersetup{colorlinks = true}

\numberwithin{equation}{section}

\allowdisplaybreaks

\def\qed{\nopagebreak\hfill{\rule{4pt}{7pt}}
\medbreak}

\newlength{\boxedparwidth}
  {\begin{center} \begin{tabular}{|@{\hspace{.315in}}c@{\hspace{.15in}}|}
                  \hline \\ \begin{minipage}[t]{\boxedparwidth}
                  \setlength{\parindent}{.25in}}%
  {\end{minipage} \\ \\ \hline \end{tabular} \end{center}}

\begin{document}
\begin{center}

 {\Large \bf Generalizations of Euler's Theorem to $k$-regular partitions}
\end{center}

\begin{center}
  {Hongshu Lin$^1$}, and {Wenston J.T. Zang$^2$} \vskip 2mm

 $^{1,2}$School of Mathematics and Statistics, Northwestern Polytechnical University, Xi'an 710072, P.R. China\\[6pt]
	$^{1,2}$ MOE Key Laboratory for Complexity Science in Aerospace, Northwestern Polytechnical University, Xi'an 710072, P.R. China\\[6pt]
	$^{1,2}$ Xi'an-Budapest Joint Research Center for Combinatorics, Northwestern Polytechnical University, Xi'an 710072, P.R. China\\[6pt]

   \vskip 2mm

$^1$linhongshu@mail.nwpu.edu.cn, 
    $^2$zang@nwpu.edu.cn
\end{center}

\vskip 6mm \noindent {\bf Abstract.} Let $A_k(n)$ denote the set of $k$-distinct partitions of $n$, and let $B_k(n)$ be the set of $k$-regular partitions of $n$. Glaisher showed that $\# A_k(n) = \# B_k(n)$. For $k=2$, this equality yields the celebrated Euler's partition theorem. In this paper, we present a new partition set $E_k(n)$, which is equinumerous to $B_k(n)$.



\section{Introduction}

The origin of all partition identities is Euler's theorem\cite{Andrews-1998}, and we know that
\begin{equation}\label{equ-ab}
    \#A(n)=\#B(n),
\end{equation}
where $A(n)$ is the set of partitions of $n$ with distinct parts and $B(n)$ is the set of partitions of $n$ with each part odd.

% Furthermore, as the extension on Euler's Theorem, Glaisher\cite{Glaisher-1883} gave that 
% \begin{equation}
%     \#A_k(n)=\#B_k(n),
% \end{equation}
% where $A_k(n)$ is the set of partitions of $n$ with each part appearing not exceeding $k-1$ times and $B_k(n)$ is the set of partitions of $n$ with each part is not the multiples of $k$.

Euler's theorem encompasses all the key elements conducive to generalization, and over the centuries, relevant generalizations have emerged in abundance. The Rogers-Ramanujan identities\cite{Ramanujan-1927,Ramanujan-1988,Rogers-1894} and Schur’s theorem\cite{Schur-1926} laid the foundation for twentieth-century advancements in this field; for a comprehensive overview of some results from the late twentieth century, refer to Henry Alder’s survey\cite{Alder-1969}.

Recently, Andrews, Kumar and Yee\cite{Andrews-2025} have given that the extension of Euler's Theorem itself with two further partition functions. 
\begin{thm}[\cite{Andrews-2025}]\label{thm-00}
     For $n\geq 0$, we have 
    \begin{equation}
        \#A(n)=\#B(n)=\#C(n+1)=\frac{1}{2}\#D(n+1),
    \end{equation}
    where $C(n)$ is the set of partitions of $n$ with largest part even and parts not exceeding half of the largest part are distinct, and $D(n)$ is the set of partitions of $n$ into non-negative parts wherein the smallest part appear exactly twice and no other parts are repeated.
\end{thm}

The main purpose of this paper is to  generalize Theorem \ref{thm-00} to $k$-regular partitions. Recall that $k$-regular partitions and $k$-distinct partitions were first introduced by Glaisher\cite{Glaisher-1883} as follows:

\begin{thm}[\cite{Glaisher-1883}]\label{thm-gla}
    For $n\geq 0$, let $A_k(n)$ denote the set of $k$-distinct partitions of $n$, that is the appearances of each number  not exceeding $k-1$ times, and let $B_k(n)$ denote the set of $k$-regular partitions of $n$, which means each part is not the multiples of $k$. Then we have
    \begin{equation}\label{equ-akbk}
        \#A_k(n)=\#B_k(n).
    \end{equation}
\end{thm}

In recent years, researchers have conducted extensive studies on $k$-distinct partitions and $k$-regular partitions. Arithmetic properties of $k$-regular partition functions have received a great deal of attention (see, for example \cite{Andrews-2010,Gordon-1997,Hirschhorn-2010,Ono-2000}); Many scholars have conducted extensive research on  inequalities of hook length on $k$-regular partitions and $k$-distinct partitions (see, for example \cite{Ballantine-2023,Kim,Singh-Barman-2024}).

Let $E_k(n)$ denote the set of partitions of $n$ with the largest part $m$ is not  multiples of $k$ and parts not exceeding $\lfloor\frac{m}{k}\rfloor-1$ appear at most $k-1$ times. The main result of this paper can be stated as follows.

\begin{thm}\label{thm-4}
    For $n\geq 1$, we have
    \begin{equation}
       \#B_k(n)=\#E_k(n).
    \end{equation}
\end{thm}

Setting $k=2$ in Theorem \ref{thm-4},  we obtain a new type of partition which is equinumerous to $B(n)$.
\begin{core}\label{thm-3}
    For $n\geq 0$, we have
    \begin{equation}
        \#B(n)=\#E(n).
    \end{equation}
    where $E(n)$ is the set of partitions of $n$ with largest part $m$ odd and parts not exceeding $\frac{m-1}{2}$ are distinct.
\end{core}

In addition, we also introduce a partition function $C_k(n)$, which is a $k$-distinct generalization of $C(n)$. Let $C_k(n)$ be the set of the partition of $n$ with largest part $m\equiv 0 \pmod{k}$ and parts not exceeding $\frac{m}{k}$ appear at most $k-1$ times, and let  $B'_k(n)$ be the set of the $k$-regular partition of $n$ with the largest part $m\equiv-1 \pmod{k}$. We have the following result. 
\begin{thm}\label{thm-2}
    For $n\geq 0$, we have 
    \begin{equation}
        \#B'_k(n)=\#C_k(n+1).
    \end{equation}
\end{thm}

Clearly set $k=2$, we recover $\#B(n)=\#C(n+1)$.

This paper is organized as follows. We give the combinatorial proof of  Theorem \ref{thm-4} in Section \ref{sec2}. Section \ref{sec3} is devoted to giving a generating function proof and a combinatorial proof of Theorem \ref{thm-2}.

\section{Proof of Theorem \ref{thm-4}}\label{sec2}

In this section, we build a bijection between $B_k(n)$ and $E_k(n)$, which proves Theorem \ref{thm-4}.

We first introduce some notations of partitions. Let $f_\lambda(t)$ be the number of appearances of $t$ in $\lambda$, and we write $\lambda$=$(\lambda_1,\lambda_2,\ldots,\lambda_\ell)$ as $(1^{f_\lambda(1)},2^{f_\lambda(2)},3^{f_\lambda(3)},\ldots)$. When $f_\lambda(i)=1$, we write $i$ instead of $i^1$ for simplicity. Moreover, when $f_\lambda(i)=0$, we may omit the term $i^0$. Furthermore, we let $\ell(\lambda)$ denote the number of parts of $\lambda$. For example, $\lambda=(9,3,3,2,1,1,1)=(1^3,2,3^2,9)$ and $\ell(\lambda)=7$.

Let $\lambda,\mu$ be two partitions. We use $\lambda\cup\mu$ to denote the partition
\[(1^{f_\lambda(1)+f_\mu(1)},2^{f_\lambda(2)+f_\mu(2)}\ldots,(m^{f_\lambda(m)+f_\mu(m)}),\]where $m=\max\{\lambda_1,\mu_1\}$. For example, $\lambda=(1^5,2^4,5^2,6)$ and $\mu=(1^2,2^3,5)$, we have $\lambda\cup\mu=(1^7,2^7,5^3,6)$.

Now we give the bijection $\psi$ between $B_k(n)$ and $E_k(n)$, which yields a combinatorial proof of Theorem \ref{thm-4}.

{\noindent\it Proof of Theorem \ref{thm-4}.} Given the partition $\lambda=(\lambda_1,\lambda_2,\ldots,\lambda_\ell)\in E_k(n)$, then by definition $\lambda$ is a partition of $n$ with $k\nmid\lambda_1$. Moreover, assume $\lambda_1=ki-r$, where $i\ge 1$ and $1\le r\le k-1$, then each part not exceeding $i-1$ appears at most $k-1$ times.  We first define a map $\alpha$ from $\mathbb{N}$ to the set of partitions as follows. For any $t\ge 1$, write $t=k^a\cdot b$, where $k\nmid b$ and $a\ge 0$. Then define
\begin{equation}
    \alpha(t)=\underset{k^a\text{'s }b}{(\underbrace{b,b,\ldots,b})}.
\end{equation}
Clearly $t=|\alpha(t)|$.

We next give the bijection $\psi$ between $E_k(n)$ and $B_k(n)$. Define
\begin{equation}
    \mu=\psi(\lambda)=\bigcup_{j=1}^\ell\alpha(\lambda_j).
\end{equation}

Since $k\nmid b$, we see that each part in $\alpha(t)$ is not a multiple of $k$. Furthermore, $$|\mu|=\sum_{j=1}^\ell |\alpha(\lambda_j)|=\sum_{j=1}^\ell|\lambda_j|=n.$$
 Thus $\mu\in B_k(n)$. 

We next show that $\psi$ is a bijection. To this end, we establish a map $\psi^{-1}$ from $B_k(n)$ to $E_k(n)$. Given $\mu\in B_k(n)$, assume $\mu=(1^{f_{\mu}(1)},2^{f_{\mu}(2)},\ldots,n^{f_{\mu}(n)})$. By definition we see that $f_\mu(t)=0$ whenever $k\mid t$.  For the largest part $\mu_1$, we set that $\mu_1=ki-r$ where $i$ is a positive integer and $1\leq r\leq k-1$. For each $s^{f_\mu(s)}$, where $k\nmid s$, we define a map $\beta$ that transforms the block $(s^{f_\mu(s)})$ into a new partition as follows.

 For each $s$, write $f_\mu(s)$ in base $k$-form, that is 
 \begin{equation}\label{equ-k-base}
     f_\mu(s)=f_0+f_1k+f_2k^2+\cdots+f_mk^m,
 \end{equation} 
 where $0\leq f_j \leq k-1$ for each $0\le j\le m$. Let $t=\min\{j:k^js\geq i,j\geq 0\}$. Then define 
    $$\beta(s^{f_\mu(s)})=(s^{f_0},(ks)^{f_1},(k^2s)^{f_2},\ldots,(k^{t-1}s)^{f_{t-1}},(k^ts)^{\sum_{j=0}^{m-t}k^j f_{t+j}}).$$

It is clear that
\begin{align}\label{eq-beta-weight}
    |\beta(s^{f_\mu(s)})|&=f_0s+f_1ks+f_2k^2s+\cdots+f_{t-1}k^{t-1}s+\left(\sum_{j=0}^{m-t}k^jf_{t+j}\right)k^ts\nonumber\\
    &=\left(f_0+f_1k+\cdots+f_{t-1}k^{t-1}+\sum_{j=t}^{m}k^jf_{j}\right)s=f_\mu(s)s.
\end{align}

We next give the bijection $\psi^{-1}$ between $B_k(n)$ and $E_k(n)$. Define
\begin{equation}
    \lambda=\psi^{-1}(\mu)=\bigcup_{s=1}^n\beta(s^{f_\mu(s)}).
\end{equation}

It is obvious that for any $s\ne s_1$, $k\nmid s$ and $k\nmid s_1$, we have $k^js\ne k^{j_1}s_1$ for any $j$ and $j_1$. Thus  the part $k^j s$ is produced exclusively by $\beta(s^{f_\mu(s)})$. Therefore when $k^js\le i-1$, we see that $f_\lambda(k^js)=f_j\le k-1$, where $f_j$ is defined in \eqref{equ-k-base}. 

Moreover, from the choice of $t$, we see that either both $t=0$ and $\mu_1\ge s\ge i$ or $ k^{t-1}s \le i-1$ which means $k^ts\le ki-k<\mu_1$. So in either case we have $k^ts\le \mu_1$, which implies $\lambda_1=\mu_1$. Thus $k\nmid \lambda_1=\mu_1=ki-r$. Moreover, by \eqref{eq-beta-weight} we see that $|\lambda|=|\mu|$. This implies that $\mu\in E_k(n)$. Furthermore, verification shows that $\psi^{-1}(\psi(\lambda))=\lambda$ and $\psi(\psi^{-1}(\mu))=\mu$, which proves $\psi$ is a bijection between $B_k(n)$ and $E_k(n)$.\qed

For example, let $k=3$ and $\lambda=(7,6,6,3,2,2,1,1)\in E_3(28)$, using the bijection $\psi$, we have $\mu=\psi(\lambda)=(7,2^8,1^5)\in B_3(28)$. Applying $\psi^{-1}$ on $\mu$, we recover $\lambda$.

\section{Proof of Theorem \ref{thm-2}}\label{sec3}

In this section, we provide two proofs of Theorem \ref{thm-2}. We begin with a proof via generating functions, followed by a combinatorial proof.

{\noindent \it Proof of Theorem \ref{thm-2}.} We first give the generating function of $B'_k(n)$ and $C_k(n)$ as follows.

  Let $B'_{k,i}(n)$ denote the subset of $B'_k(n)$ with the largest part equal to $ki-1$. Its generating function is given by
\begin{equation}
    \sum_{n=0}^\infty B'_{k,i}(n)q^n=\frac{(q^{k};q^k)_{i-1}q^{ki-1}}{(q;q)_{ki-1}}.
\end{equation}

Here we use the standard $q$-series notation
\[
(a;q)_n=\prod_{i=1}^{n}(1-aq^{i-1})\quad \text{and} \quad
(a;q)_\infty=\prod_{i=1}^{\infty}(1-aq^{i-1}).
\]
Summing over all $i \ge 1$, we obtain the generating function for $B'_k(n)$:
    \begin{align*}
    \sum_{n=0}^\infty B'_k(n)q^n=\sum_{i=1}^\infty\sum_{n=0}^\infty B'_{k,i}(n)q^n
    =\sum_{i=1}^\infty\frac{(q^{k};q^k)_{i-1}q^{ki-1}}{(q;q)_{ki-1}}.
    \end{align*}   

 Similarly, let $C_{k,i}(n)$ denote the subset of $C_k(n)$ with largest part $ki$, then by definition each number not exceeding $i$ appear at most $k-1$ times. Thus 
\begin{equation}
    \sum_{n=0}^\infty C_{k,i}(n)q^n=\prod_{j=1}^i\left(\sum_{l=0}^{k-1}q^{jl}\right)\frac{q^{ki}}{(q^{i+1};q)_{(k-1)i}}=\frac{(q^k;q^k)_i}{(q;q)_i}\frac{q^{ki}}{(q^{i+1};q)_{(k-1)i}}.
\end{equation}
Summing over $i\ge 0$ yields 
\begin{align*}
    \sum_{n=0}^\infty C_k(n)q^n=\sum_{i=0}^\infty\sum_{n=0}^\infty C_{k,i}(n)q^n
    =\sum_{i=0}^\infty \frac{(q^k;q^k)_i}{(q;q)_i}\frac{q^{ki}}{(q^{i+1};q)_{(k-1)i}}.
\end{align*}
We now relate the two generating functions. Starting from the expression for $B'_k(n)$:
\begin{align*}
    \sum_{n=1}^\infty B'_k(n)q^{n+1}
    =&\sum_{i=1}^\infty \frac{(q^{k};q^k)_{i-1}q^{ki}}{(q;q)_{ki-1}}\\
    =&\sum_{i=1}^\infty \frac{(q^{k};q^k)_{i}q^{ki}}{(q;q)_{ki}}\\
    =&\sum_{i=0}^\infty \frac{(q^{k};q^k)_{i}q^{ki}}{(q;q)_{ki}}-1\\
    =&\sum_{i=0}^\infty \frac{(q^k;q^k)_iq^{ki}}{(q;q)_i(q^{i+1};q)_{(k-1)i}}-1=\sum_{n=0}^\infty C_k(n)q^n-1.
\end{align*}
   Comparing the coefficients of $q^{n+1}$, we have   $B'_k(n)=C_k(n+1)$.\qed

We conclude this paper by   presenting  a combinatorial proof of Theorem \ref{thm-2}.

{\noindent \it Combinatorial proof of Theorem \ref{thm-2}.} Let $\lambda \in C_k(n+1)$. By definition,  $\lambda_1=ki$ for some $i\ge 1$ and  $f_\lambda(j)\le k-1$ for each $1\le j\le i$. Now, we define a map $\alpha$ from $\mathbb{N}$ to the set of partitions as follows. 

For the largest part $\lambda_1=ki$, define $\alpha(\lambda_1)=(ki-1)$. For $\ell(\lambda)\ge j\ge 2$, write  $\lambda_j=k^{a}\cdot b$ where $k\nmid b$, then define
$$\alpha(\lambda_j)=\underset{k^a\text{'s } b}{(\underbrace{b,b,\ldots,b})}.$$

We next give the bijection $\phi$ between $C_{k}(n+1)$ and $B'_k(n)$. Define
\begin{equation}
    \mu=\phi(\lambda)=\bigcup_{j=1}^{\ell(\lambda)}\alpha(\lambda_j).
\end{equation}

We first show that $\mu \in B'_k(n)$. By construction of $\phi$, the largest part of $\mu$ is $\mu_1=ki-1 \equiv -1 \pmod{k}$ in $\mu$. Moreover, $\alpha(\lambda_j)$ is a partition with each part not divisible by $k$, and so is $\mu$. Furthermore, $$|\mu|=\sum_{j=1}^\ell|\lambda_j|-1=n+1-1=n.$$
Thus $\mu \in B'_k(n)$.

To show that $\phi$ is a bijection, we construct its inverse $\phi^{-1}: B'_k(n) \to C_k(n+1)$. Given $\mu \in B'_k(n)$, by definition $\mu_1=ki-1$ for some $i\ge 1$. Moreover,  $k\nmid\mu_j$ for any $1\le j\le \ell(\mu).$ For each $k\nmid s$, we define a map $\beta$ from $s^{f_\mu(s)}$ to the set of partitions as follows.

For the largest part $\mu_1=ki-1$, define that $ \beta(\mu_1)=(ki)$. Set $\nu=(\mu_2,\ldots,\mu_{\ell(\mu)})$, for each part $s \in \nu$, write 
$$f_\nu(s)=f_0+f_1k+f_2k^2+\cdots +f_mk^m,$$ 
where $0\leq f_j\leq k-1$ for each $0\leq j\leq m$. Assume $t=\min\{j:k^js >i,j\geq 0\}$, set
$$\beta(s^{f_\nu(\nu_j)})=(s^{f_0},(ks)^{f_1},(k^2s)^{f_2},\ldots,(k^{t-1}s)^{f_{t-1}},(k^ts)^{\sum_{j=0}^{m-t} k^jf_{t+j}}).$$

It is clear that
\begin{align*}
    |\beta(s^{f_\nu(s)})|&=f_0s+f_1ks+f_2k^2s+\cdots+f_{t-1}k^{t-1}s+\left(\sum_{j=0}^{m-t}k^jf_{t+j}\right)k^ts\\
    &=\left(f_0+f_1k+\cdots+f_{t-1}k^{t-1}+\sum_{j=t}^{m}k^jf_{j}\right)s=f_\nu(s)s.
\end{align*}

We next give the bijection $\phi^{-1}$ between $B'_k(n)$ and $C_k(n+1)$. Define
\begin{equation}
    \lambda=\phi^{-1}(\mu)=\bigcup_{s=1}^n\beta(s^{f_\nu(s)})\cup\beta(\mu_1).
\end{equation}

Using the same argument as in the proof of Theorem~\ref{thm-4}, one may verify that for any $\lambda_j \le i$, we have $f_\lambda(\lambda_j) \le k-1$. Moreover, with  $\beta(\mu_1)=ki$ divided by $k$ and $i<k^t\nu_j\leq ki$, we see that  $\lambda_1=ki$. Furthermore
$$|\lambda|=\sum_{\nu_j\in\nu}\beta(\nu_j^{f_\nu(\nu_j)})+\beta(\mu_1)=|\nu|+\mu_1+1=n+1. $$
Hence $\mu\in C_k(n+1)$. Furthermore, we can check that $\phi^{-1}(\phi(\lambda))=\lambda$ and $\phi^{-1}(\phi(\mu))=\mu$. Thus $\phi$ is a bijection between $B'_k(n)$ and $C_k(n+1)$.\qed

For example, let $k=3$ and $\lambda=(9,9,3,2,2,1,1) \in C_3(27)$, using the bijection $\phi$, we have that $\mu=\phi(\lambda)=(8,2,2,1^{14})\mu\in B'_3(26)$. Applying $\phi^{-1}$ on $\mu$, we recover $\lambda$.

\noindent{\bf Acknowledgments.}   This work was supported by the National Science Foundation of China grants 12171358 and 12371336.

\end{document}